\documentclass[12pt,final]{amsart}

\usepackage{verbatim}

\usepackage{amssymb}
\usepackage{amscd}

\newtheorem{lemma}{Lemma}[section]
\newtheorem{proposition}[lemma]{Proposition}
\newtheorem{theorem}[lemma]{Theorem}

\newtheorem{corollary}[lemma]{Corollary}

\theoremstyle{definition}
\newtheorem{defn}[lemma]{Definition}

\theoremstyle{remark}
\newtheorem{remark}[lemma]{Remark}
\newtheorem{example}[lemma]{Example}

\newcommand{\Z}{\ensuremath{{\mathbb Z}}}
\newcommand{\N}{\ensuremath{{\mathbb N}}}

\newcommand{\R}{\ensuremath{\mathbb R}}
\newcommand{\C}{\ensuremath{\mathbb C}}

\newcommand{\bft}{\ensuremath{{\mathbf{T}}}}
\newcommand{\bbt}{\ensuremath{{\mathbb{T}}}}

\newcommand{\hk}{\ensuremath{{\mathcal {H}}_k}}

\DeclareMathOperator{\sign}{sign}
\DeclareMathOperator{\End}{End}
\DeclareMathOperator{\td}{Td}
\DeclareMathOperator{\ch}{ch}

\author{Tatyana Barron}
\address{T. Barron, Department of Mathematics, 
University of Western Ontario, 
London, Ontario N6A 5B7, Canada }
\email{tatyana.barron@uwo.ca}

\author{Baran Serajelahi}
\address{B. Serajelahi, Department of Mathematics, 
University of Western Ontario, 
London, Ontario N6A 5B7, Canada}
\email{bserajel@uwo.ca}

\thanks{Research of the first author is supported in part 
by NSERC}

\title[Berezin-Toeplitz quantization, hyperk\"ahler manifolds...]{Berezin-Toeplitz quantization, hyperk\"ahler manifolds, and multisymplectic manifolds}

\begin{document}
\sloppy

\maketitle

{\bf Abstract.} We suggest a way to quantize, using Berezin-Toeplitz quantization, 
a compact hyperk\"ahler manifold (equipped with a natural $3$-plectic form), 
or a compact integral K\"ahler manifold of complex dimension $n$ regarded as a $(2n-1)$-plectic manifold. We show that quantization has reasonable semiclassical properties.

\section{Introduction}

(Berezin-)Toeplitz quantization, while interesting to study by itself, 
also has turned out to be a useful tool in several areas of mathematics. 
Over the years it was found to have applications to deformation quantization 
(see e.g. \cite {schlich:00}, \cite{karabegov:01}), 
to study of the Hitchin connection and TQFT   
(work of J. Andersen, see in particular \cite{andersen:06}, \cite{andersen:10}),   
L. Polterovich's work on rigidity of Poisson brackets \cite{polterovich:12}, 
and work of Y. Rubinstein and S. Zelditch \cite{rubins:12} on homogeneous complex Monge-Amp\`ere equation, 
in connection to geodesics on the space of K\"ahler metrics.  
T. Foth (T. Barron) and A. Uribe applied Berezin-Toeplitz quantization 
to give another proof of Donaldson's "scalar curvature is a moment map" statement \cite{foth:07}.

In this paper we discuss how to use Berezin-Toeplitz quantization to quantize 
hyperk\"ahler manifolds or two types of multisymplectic manifolds. 

Geometric quantization and K\"ahler/Berezin-Toeplitz quantization associate a Hilbert space (say, 
$\mathcal{H}$)  
and operators on it to a symplectic manifold $(M,\omega)$. In physics' terminology  
this is a way to pass from classical Hamiltonian mechanics to a quantum system.
Let $C^{\infty}(M)$ denote the space of complex-valued smooth functions on $M$. 
Quantization is a linear map $C^{\infty}(M)\to \{ \text{operators on} \ {\mathcal{H}}\}$, $f\mapsto \hat{f}$,  
satisfying a version of Dirac's quantization conditions: 

$1\mapsto const(\hbar)I$, 

$\{ f,g\} \mapsto const(\hbar)[\hat{f},\hat{g}]$.

It is probably fair to say that geometric quantization was developed and mainstreamed 
in the 1950s-1960s, by representation theorists,  
including Kostant, Kirillov and others,  whose primary agenda was to look for 
representations of infinite-dimensional Lie algebras with certain  properties, 
and who found this language to be quite convenient. 

Berezin-Toeplitz quantization 
can be regarded as a version of geometric quantization. In the case 
when the symplectic manifold is, moreover, K\"ahler, 
it is also referred to as K\"ahler quantization. 
The groundwork for Berezin-Toeplitz quantization was laid in 
\cite{berezin:74}, \cite{boutet:81}.   
Well-known  
Theorem \ref{bmstheorem}(i) below shows that in the framework 
of Berezin-Toeplitz quantization  
the $\{ .,.\} \leadsto [.,.]$ quantization condition 
is satisfied in the semiclassical limit $\hbar=\frac{1}{k}\to 0$, which is essentially 
the best one can get, due to various no-go theorems.   
 
There are physical systems whose behaviour is encoded by an $m$-{\it plectic} form on $M$ 
(i.e. a closed non-degenerate $m+1$-form), $\Omega$, for $m\ge 1$. 
The case $m=1$ is when $\Omega$ is symplectic. Specific examples from physics, 
with $m\ge 2$, are discussed in \cite{nambu:73}, \cite{chat:96}, \cite{baez:10}. 
See also discussion and references in \cite{cantr:99}. 
Multisymplectic geometry has been thoroughly studied by mathematicians.  
See, in particular, \cite{martin:88}, \cite{cantr:99}, \cite{madsen:12}, \cite{bur:13},  
\cite{takhtajan:94}, \cite{baez:10}, \cite{baezrog:10}, \cite{rogers:12}. 
There has been extensive discussion of {\it quantization} of $n$-plectic manifolds 
in physics literature, and substantial amount of work has been done by mathematicians too. 
See, for example, \cite{nambu:73}, \cite{takhtajan:94}, \cite{chat:96}, \cite{dito:97},  
\cite{curtr:03}, \cite{curtr:04}, \cite{debellis:10}, \cite{samann:13}, \cite{rogers:13}, 
\cite{vaisman:99}. Work of C. Rogers \cite{rogers:13} addresses 
quantization of $2$-plectic manifolds. It seems that the appropriate quantum-mechanical setting 
there involves a category, instead of a vector space, and  
intuitively this makes sense because an (integral) $2$-plectic form corresponds 
to a gerbe and sections of a gerbe form a category, not a vector space.  

There have been attempts, informally speaking, "to embed a multisymplectic physical system  
into Hamiltonian system" \cite{bayen:75}, \cite{mukunda:76}, \cite{debellis:10}. 
As far as we know, there is no known canonical way of doing this.

DeBellis, S\"amann and Szabo \cite{debellis:10} used Berezin-Toeplitz quantization 
for multisymplectic spheres via embedding them in a certain explicit way into 
complex projective spaces ${\mathbb{CP}}^{q}$  
and using Berezin-Toeplitz quantization on ${\mathbb{CP}}^{q}$. 
This is somewhat related to our results in Section \ref{volumeform}, only for $M=S^2$ 
(because among spheres only $S^2$ admits a K\"ahler form).

Let $(M,\omega)$ be a compact connected integral K\"ahler manifold of complex dimension $n$. 
In this paper we are looking into two situations when the $m$-plectic form 
$\Omega$  on  $(M,\omega)$ is constructed from the K\"ahler form (or forms): 

(I) $m=2n-1$, $\Omega = \frac{\omega^n}{n!}$ 

(II) $M$ is, moreover, hyperk\"ahler, $m=3$, 
$$
\Omega = \omega_1\wedge \omega_1 +  \omega_2\wedge \omega_2 + \omega_3\wedge \omega_3
$$
where $\omega_1=\omega, \omega_2, \omega_3$ are the three K\"ahler forms on $M$ 
given by the hyperk\"ahler structure. 

It is well-known (and easy to prove) that a volume form on an oriented $N$-dimensional manifold 
is an $(N-1)$-plectic form, and that the $4$-form above is a $3$-plectic form 
on a hyperk\"ahler manifold. See, for example, \cite{cantr:99}, \cite{rogers:12}.  

It is intuitively clear that in these two cases the classical multisymplectic system is 
essentially built from Hamiltonian system(s) and  
it should be possible to quantize $(M,\Omega)$ using the (Berezin-Toeplitz) quantization 
of $(M,\omega)$. We discuss case (I) in section \ref{volumeform}, 
case (II) in section \ref{hyperkahler}. 
Semiclassical asymptotics are the content of Theorems \ref{thvolform}, \ref{thhyperk},  
\ref{thdim4}, \ref{tensorth}, Propositions \ref{commvolform},  \ref{commdim4}, 
\ref{tensorprop}, \ref{tensorprop4}, Corollary \ref{tensorcorcomm2}. 
In both cases there are natural multisymplectic analogues of the Poisson bracket and 
the commutator: an almost Poisson bracket $\{ .,...,.\}$ 
and the generalized commutator $[.,...,.]$. 
Our discussion mainly revolves around the $\{ .,...,.\} \leadsto [.,...,.]$ quantization condition. 

The main result of section \ref{volumeform} is Theorem \ref{thvolform}. It is an analogue, 
for brackets of order $2n$, of well-known 
Theorem \ref{bmstheorem}(i) (and of its $C^l$  analogue ($l\in \N)$ 
from \cite{barron:14}). 

In section \ref{hyperkahler} we work on a hyperk\"ahler manifold $M$.    
For a smooth function $f$ on $M$ we have three Berezin-Toeplitz 
operators 
$T_{f;1}^{(k)}$, $T_{f;2}^{(k)}$, $T_{f;3}^{(k)}$, and to four smooth functions 
$f,g,h,t$ on $M$ we associate three brackets of order $4$: $\{ f,g,h,t\} _r$, $r=1,2,3$. 
In subsection   \ref{directsumgen} we show that the direct sum of generalized 
commutators is asymptotic to 
$$
T_{\{ f,g,h,t\} _1;1 }^{(k)} \oplus T_{\{ f,g,h,t\} _2;2 }^{(k)} \oplus T_{\{ f,g,h,t\} _3;3 }^{(k)} 
$$
(Theorem \ref{thhyperk}). In subsection \ref{directsum4} we show that the attempt 
to formulate everything on {\it one} vector space (not three), by taking direct sums, 
goes through all the way in the case when $M$ is the $4$-torus with 
three linear complex structures, where 
we get a straightforward analogue of Theorem \ref{bmstheorem}(i) - see Example \ref{R4} 
(\ref{asymptorus}). 
In subsection \ref{tensorproduct} we take the tensor product of the three operators, instead.    
Tensor product of generalized commutators is asymptotic to 
$$
T_{\{ f,g,h,t\} _1;1 }^{(k)} \otimes T_{\{ f,g,h,t\} _2;2 }^{(k)} \otimes T_{\{ f,g,h,t\} _3;3 }^{(k)} 
$$
(Proposition \ref{tensorprop4}). Asymptotic properties of commutators and generalized commutators 
of operators $\bbt ^{(k)}_f= T_{f;1}^{(k)}\otimes T_{f;2}^{(k)}\otimes T_{f;3}^{(k)}$ 
are captured in Prop. \ref{tensorprop} and Theorem \ref{tensorth}.

We note that while, for simplicity, 
the exposition throughout the paper is for $C^{\infty}$ symbols, 
- all our results hold, in fact, for $C^4$ symbols. 
To modify the proofs in order to get the same statements for $C^4$  symbols, the estimates from \cite{bordemann:94} should be replaced 
by estimates from \cite{barron:14} - see subsection \ref{nonsmoothsymb}. 
Results from \cite{barron:14} allow to tackle the case of $C^2$ and $C^3$ symbols as well,  
but we do not include the corresponding version of our results 
(the asymptotics will differ from the $C^{\infty}$ case). 

This paper is a part of the Ph.D. thesis of the second author who is co-supervised 
by the first author and M. Pinsonnault. 

{\bf Acknowledgements.} We are thankful to G. Denham, M. Gualtieri, B. Hall, N. Lemire, A. Uribe, 
K. Yoshikawa, for brief related discussions, and to X. Ma and G. Marinescu - for comments.  
We are grateful to M. Pinsonnault 
for many questions and comments. 
We appreciate referee's suggestions that helped improve exposition 
in the paper. 

\section{Preliminaries} 
\subsection{Some notations and definitions} 
Throughout the paper we shall use the following notations: 

$S_n$, for a positive integer $n$, will denote the symmetric group (i.e. the group of permutations of $1,...,n$), 

for a finite-dimensional complex vector space $V$ and $A,B\in \End(V)$ 
$[A,B]=AB-BA$, 

$I$ will denote the identity operator on $V$, 

if $V$ is equipped with a norm, then $||A||$ will denote the operator norm of $A$, 

$C^{\infty}(M)$ will denote the algebra of smooth complex-valued functions 
on a smooth manifold $M$, 

for $f\in C^\infty (M)$ we write $|f|_{\infty}=\sup _{x\in M}|f(x)|$.

\begin{defn}
An $(m+1)$-form $\Omega$ on a smooth manifold  $M$ is called an {\bf m-plectic form} 
if it is closed (i.e. $d\Omega =0$) and non-degenerate (i.e. 
$v\in T_xM, v\lrcorner \Omega _x=0 \Rightarrow v=0$). 

If $\Omega$ is an $m$-plectic form on $M$, $(M,\Omega )$ is called a {\bf multisymplectic}, 
or {\bf m-plectic}, manifold.  
\end{defn}

\begin{defn}(\cite{takhtajan:94}, \cite{gautheron:96}) 
Let $M$ be a smooth manifold.  
A multilinear map 
$$
\{ .,...,.\} :(C^{\infty}(M))^{\otimes j}\to C^{\infty}(M)
$$ 
is called a {\bf Nambu-Poisson bracket} or {\bf (generalized) Nambu bracket 
of order $j$} if it satisfies the following properties: 
\begin{itemize}
\item (skew-symmetry) $\{ f_1,...,f_j\} 
= \sign (\sigma )\{ f_{\sigma (1)},...,f_{\sigma (j)}\}$ 
for all $f_1,...,f_j\in C^{\infty}(M)$ and for all $\sigma\in S_j$, 
\item 
(Leibniz rule) 
$\{ f_1,...,f_{j-1}, g_1g_2\} = \{ f_1,...,f_{j-1}, g_1\} g_2+g_1 \{ f_1,...,f_{j-1},g_2\}$
for all $f_1,...,f_{j-1}, g_1,g_2\in C^{\infty}(M)$, 
\item (Fundamental Identity) 
$$
\{ f_1,...,f_{j-1,},\{ g_1,..., g_j\} \} = 
\sum _{i=1}^j\{ g_1,...,\{ f_1,...,f_{j-1}, g_i\} ,...,g_j\} ,
$$
for all $f_1,...,f_{j-1}, g_1,...,g_j \in C^{\infty}(M)$. 
\end{itemize}
\end{defn}
It is natural to ask how to generalize the Hamiltonian formalism of symplectic geometry 
to the multisymplectic setting. We do not need the full multisymplectic formalism for the purposes 
of this paper, and we refer the reader to \cite{takhtajan:94}, \cite{helein:04}, \cite{rogers:12}. 

\begin{defn}\label{defgpb}
(\cite{azcar:96}, \cite{azcar:10}) 
Let $M$ be a smooth manifold and suppose $j$ is an even positive integer.  
A multilinear map 
$$
\{ .,...,.\} :(C^{\infty}(M))^{\otimes j}\to C^{\infty}(M)
$$ 
is called a {\bf generalized Poisson bracket}  
if it satisfies the following properties: 
\begin{itemize}
\item (skew-symmetry) $\{ f_1,...,f_j\} 
= \sign (\sigma )\{ f_{\sigma (1)},...,f_{\sigma (j)}\}$ 
for all $f_1,...,f_j\in C^{\infty}(M)$ and for all $\sigma\in S_j$, 
\item 
(Leibniz rule) 
$\{ f_1,...,f_{j-1}, g_1g_2\} = \{ f_1,...,f_{j-1}, g_1\} g_2+g_1 \{ f_1,...,f_{j-1},g_2\}$
for all $f_1,...,f_{j-1}, g_1,g_2\in C^{\infty}(M)$, 
\item (Generalized Jacobi Identity) 
$$
\mbox{Alt}\{ f_1,...,f_{j-1},\{ f_j,..., f_{2j-1}\} \} = 
$$
$$
\sum _{\sigma\in S_{2j-1}} \sign(\sigma) \{ f_{\sigma(1)},...,f_{\sigma(j-1)},\{ 
f_{\sigma(j)},..., f_{\sigma(2j-1)}\} \} =0
$$
for all $f_1,...,f_{2j-1} \in C^{\infty}(M)$. 
\end{itemize}
\end{defn}
\begin{defn} (\cite{ibanez:97})
A bracket as in Definition \ref{defgpb} satisfying only the first two conditions (skew-symmetry and Leibniz rule) 
is called an {\bf almost Poisson bracket of order j}. 
\end{defn}
\begin{remark}
A Nambu-Poisson bracket of even order is a generalized Poisson bracket \cite{ibanez:97}. 
\end{remark}
\subsection{Generalized commutator} 
Let $[.,.,.,.]$ denote the Nambu generalized commutator (\cite{nambu:73}, \cite{takhtajan:94}, \cite{chat:96}):
for a finite-dimensional complex vector space $V$ and $A_1,...,A_{2n}\in \End(V)$
$$
[A_1,...,A_{2n}]=\sum_{\sigma\in S_{2n}} \sign(\sigma )A_{\sigma(1)}...A_{\sigma(2n)} . 
$$
For example, for $n=2$ 
$$
[A_1,A_2,A_3,A_4]=\sum_{\sigma\in S_4} \sign(\sigma)A_{\sigma(1)}A_{\sigma(2)}A_{\sigma(3)}
A_{\sigma(4)} =
$$
\begin{equation}
\label{comm4}
\begin{split}
[A_1,A_2][A_3,A_4] - [A_1,A_3][A_2,A_4] + [A_1,A_4][A_2,A_3]+
\\
[A_3,A_4][A_1,A_2] - [A_2,A_4][A_1,A_3] + [A_2,A_3][A_1,A_4] . 
\end{split}
\end{equation}
The bracket $[.,.,.,.]$ defines a map $\bigwedge^4 \End(V)\to \End(V)$ which does not satisfy
the Leibniz rule and does not satisfy the Fundamental Identity.
There has been some discussion of this in physics literature (e.g. \cite{curtr:03}) 
and they seem to think
that requiring these two conditions is not necessary. There has been investigation
into algebraic properties of this bracket - see e.g. \cite{curtr:09} and \cite{azcar:10}, where some ideas go back to \cite{bremner:98}, \cite{filippov:85}, and earlier work by Kurosh 
and his school.  

Let us denote, for convenience, 
$$
\sideset{}{'} 
\sum _{\sigma\in S_{2n}} =\sum_{\substack{{\sigma\in S_{2n}}\\ 
{ \sigma (1)<\sigma(2),...,} \\ {\sigma(2n-1)<\sigma(2n)}} } . 
$$
\begin{lemma} 
\label{lemcomm1}
$$
[A_1,...,A_{2n}]=\sideset{}{'}
\sum _{\sigma\in S_{2n}} \sign (\sigma ) 
[A_{\sigma(1)},A_{\sigma(2)}][A_{\sigma(3)},A_{\sigma(4)}]... 
[A_{\sigma(2n-1)},A_{\sigma(2n)}] .
$$
\end{lemma}
\noindent {\bf Proof.} Each monomial from the left hand side appears in the right hand side, 
exactly once, with the same sign. Each term from the right hand side appears 
in the left hand side. Therefore the expressions are identical.    
$\Box$ 
\begin{lemma} 
\label{lemcomm2}
$$
[A_1,...,A_{2n}]=\frac{1}{2^n}
\sum _{\sigma\in S_{2n}} \sign (\sigma ) 
[A_{\sigma(1)},A_{\sigma(2)}][A_{\sigma(3)},A_{\sigma(4)}]... 
[A_{\sigma(2n-1)},A_{\sigma(2n)}] .
$$
\end{lemma}
\noindent {\bf Proof.} By straightforward comparison of the polynomials. 
Observe that each monomial from the left-hand side appears 
in the sum in the right-hand side exactly $2^n$ times, with appropriate sign, 
and this accounts for all the terms in the right hand side. $\Box$ 
\begin{remark}
Equality (\ref{comm4}) is (93) \cite{curtr:03}. It is not hard to see that 
Lemma \ref{lemcomm2} is equivalent to (94) \cite{curtr:03}. 
\end{remark}
\subsection{Berezin-Toeplitz operators} 
Suppose $(M,\omega )$ is a compact connected K\"ahler manifold 
and the K\"ahler form $\frac{\omega}{2\pi}$ is integral. Let $L$ be a holomorphic hermitian line bundle 
such that the curvature of the Chern connection is $-i\omega$. 
Let $k$ be a positive integer. The space $H^0(M,L^{\otimes k})$ of holomorphic sections 
of $L^{\otimes k}$ is a finite-dimensional complex vector space. 
Let $\Pi_k$ denote the orthogonal projection from $L^2(M,L^{\otimes k})$ onto 
$H^0(M,L^{\otimes k})$ (the Hermitian 
inner product is obtained from the hermitian metric on $L$). 
\subsubsection{Smooth symbol}
\label{smoothsymb}

Reference used throughout this subsection is \cite{bordemann:94}, 
where the method is based on the analysis of Toeplitz structures 
from \cite{boutet:81}. Results mentioned here and more extensive 
discussion can be found in surveys on Berezin-Toeplitz 
quantization, -  for example in \cite{schlich:10}. 

For $f\in C^{\infty}(M)$ the operator 
$$
T_f^{(k)}=\Pi_k\circ (mult. \ by \ f)\in \End ( H^0(M,L^{\otimes k}) ) , 
$$
or the operator $\oplus T_f^{(k)}$,  
is called the {\it Berezin-Toeplitz operator for $f$}. 
Here are some properties of these operators that will be most frequently 
used in this paper.   

For $\alpha,\beta\in \C$ and $f,g\in C^{\infty}(M)$  
$$
T_{\alpha f+\beta g}^{(k)}=\alpha T_f^{(k)}+\beta T_g^{(k)} .
$$ 
\begin{theorem}[\cite{bordemann:94} Th. 4.1, 4.2; \cite{ma:07}, 
\cite{ma:08}] 
\label{bmstheorem}

For $f,g\in C^{\infty}(M)$, as $k\to\infty$,
\begin{itemize}
\item[(i)]
$$
||ik[T_f^{(k)},T_g^{(k)}]-T_{\{ f,g\} }^{(k)}||=O(\frac{1}{k}), 
$$
\item[(ii)]
there is a constant $C=C(f)>0$ such that 
$$
|f|_{\infty}-\frac{C}{k}\le ||T_f^{(k)}||\le |f|_{\infty} . 
$$ 
\end{itemize}
\end{theorem}
\begin{proposition}[\cite{bordemann:94} p. 291] 
\label{bmsprop}
For $f_1,...,f_p\in C^{\infty} (M)$ 
$$
||T_{f_1}^{(k)}...T_{f_p}^{(k)}-T_{f_1...f_p}^{(k)}||=O(\frac{1}{k}) 
$$
as $k\to\infty$. 
\end{proposition}
\begin{proposition}[\cite{bordemann:94} p. 289]
\label{propcomm}

For $f,g\in C^{\infty}(M)$
$$
\lim_{k\to \infty} ||[T_f^{(k)},T_g^{(k)}]||=0. 
$$
\end{proposition}
\begin{remark}
\label{remarkcomm}
Proof of this Proposition (it's one line, use Theorem \ref{bmstheorem} 
and triangle inequality) actually implies that 
$$
||[T_f^{(k)},T_g^{(k)}]||=O(\frac{1}{k})
$$
as $k\to\infty$. 
\end{remark}
\subsubsection{$C^l$ symbol}
\label{nonsmoothsymb}

The reference for theorems analogous to those above in subsection \ref{smoothsymb}, 
with $f\in C^l(M)$, is \cite{barron:14}. In \cite{barron:14} the method is different from  
\cite{bordemann:94}. It relies on techniques developed in \cite{ma:07}, \cite{ma:08}, 
see also \cite{ma:11}.  For $l=4$ statements similar to Theorem \ref{bmstheorem}, 
Prop. \ref{bmsprop} follow from Cor. 4.5, Remark 5.7(b), Cor. 4.4 of \cite{barron:14}. 
The fact that for $f,g\in C^4(M)$ $||[T_f^{(k)},T_g^{(k)}]||=O(\frac{1}{k})$
as $k\to\infty$ easily follows too, from Cor. 4.5 and  Remark 5.7(b) \cite{barron:14}.

\section{Quantization of the $(2n-1)$-plectic structure on an $n$-dimensional K\"ahler manifold}
\label{volumeform} 

Let $(M,\omega)$ be a compact connected $n$-dimensional K\"ahler manifold ($n\ge 1$). 
We shall denote by $\{ .,.\}$ the Poisson bracket for $\omega$. 
Assume that the K\"ahler form $\frac{\omega}{2\pi}$ is integral. Let $L$ be a hermitian holomorphic line 
bundle on $M$ such that the curvature of the Chern connection is equal to $-i \omega$. 

It is clear that the volume form 
$\Omega = \frac{\omega^n}{n!}$ is a $(2n-1)$-plectic form. The bracket 
$\{ .,...,.\} :\bigwedge ^{2n}C^{\infty}(M)\to C^{\infty}(M)$ defined by 
$$
df_1\wedge ...\wedge df_{2n}=\{ f_1,...,f_{2n}\} \Omega
$$
is a Nambu-Poisson bracket \cite[Cor. 1 p. 106]{gautheron:96} . 

\begin{lemma} 
\label{lembrackets}
For $f_1,...,f_{2n}\in C^{\infty}(M)$  
\begin{equation}
\label{eqfunbrackets}
\{ f_1,...,f_{2n}\} = \frac{1}{2^n n!}\sum_{\sigma\in S_{2n}}\sign (\sigma)
\prod_{j=1}^n \{ f_{\sigma(2j-1)}, f_{\sigma(2j)}\}
\end{equation}
\end{lemma}
\begin{remark}
In particular, for $n=2$ 
$$
\{f_1,f_2,f_3,f_4\} =\{ f_1,f_2\} \{ f_3,f_4\} - \{ f_1,f_3\} \{ f_2,f_4\} +  
\{ f_1,f_4\} \{ f_2,f_3\}  . 
$$
\end{remark}
\begin{remark}
For $M=\R ^{2n}$ with the standard symplectic form equality (\ref{eqfunbrackets}) 
is (7) in \cite{curtr:03}. 
\end{remark}
\noindent {\bf Proof of Lemma \ref{lembrackets}.} 
Let's use Darboux theorem and compare the left-hand side and the right-hand side 
of (\ref{eqfunbrackets}) in a local chart with coordinates 
$x_1$,...,$x_{2n}$ such that in this chart $\omega = \sum_{j=1}^n dx_{2j-1}\wedge dx_{2j}$. 
Locally, in this chart, the Poisson bracket of $f_i$, $f_l$, for $i,l\in \{ 1,...,2n\}$, is 
$$
\{ f_i,f_l\} = \sum _{j=1}^n(\frac{\partial f_i}{\partial x_{2j-1}} \frac{\partial f_l}{\partial x_{2j}}-
\frac{\partial f_i}{\partial x_{2j}} \frac{\partial f_l}{\partial x_{2j-1}} )  
$$
and 
$\{ f_1,...,f_{2n}\} =\det J$, where $J=(\frac{\partial f_i}{\partial x_l})$. 
$\det$ is the only function on $(2n)\times (2n)$ complex matrices 
which takes value $1$ on the identity matrix, linear in the rows, and takes value 
zero on a matrix whose two adjacent rows are equal  
(axiomatic characterization of the determinant, see e.g. Theorem 1.3.(3.14) \cite{artin:91}). 
The right-hand side of (\ref{eqfunbrackets}) is a polynomial  
in the entries of $J$ that 
satisfies these three conditions, therefore it must be equal to $\det J$. 
$\Box$

The following theorem shows that, informally speaking, 
$\{ .,...,.\}\to [.,...,.]$ as $k\to\infty$.  
\begin{theorem}
\label{thvolform}
For $f_1,...,f_{2n}\in C^{\infty} (M)$
$$
||\frac{(ik)^n}{n!}[T_{f_1}^{(k)},...,T_{f_{2n}}^{(k)}]-T_{\{ f_1,...,f_{2n}\} }^{(k)}|| = 
O(\frac{1}{k})
$$
as $k\to\infty$.
\end{theorem}
\noindent {\bf Proof.} 
By Theorem \ref{bmstheorem} (i) 
\begin{equation}
\label{asympcomm}
||ik[T_{f_{2j-1}}^{(k)},T_{f_{2j}}^{(k)}]-T_{\{ f_{2j-1},f_{2j}\} }^{(k)}||=O(\frac{1}{k})
\end{equation} 
for $j=1,...,n$. 
Using Prop. \ref{bmsprop} and the triangle inequality, we get: 
$$
||(ik)^n[T_{f_1}^{(k)},T_{f_2}^{(k)}]...[T_{f_{2n-1}}^{(k)},T_{f_{2n}}^{(k)}]
-T_{\{ f_1,f_2\} ...\{ f_{2n-1},f_{2n}\} }^{(k)}|| \le 
$$
$$ 
||(ik)^n[T_{f_1}^{(k)},T_{f_2}^{(k)}]...[T_{f_{2n-1}}^{(k)},T_{f_{2n}}^{(k)}]
-T_{\{ f_1,f_2\} } ^{(k)}...T_ {\{ f_{2n-1},f_{2n}\} }^{(k)}||+
$$
$$
|| T_{\{ f_1,f_2\} ...\{ f_{2n-1},f_{2n}\} }^{(k)}-
T_{\{ f_1,f_2\} } ^{(k)}...T_ {\{ f_{2n-1},f_{2n}\} }^{(k)}|| = 
$$
$$
||\bigl ( (ik[T_{f_1}^{(k)},T_{f_2}^{(k)}]-T_{\{ f_1,f_2\} } ^{(k)})+T_{\{ f_1,f_2\} } ^{(k)}\bigr ) ...
\bigl ( (ik[T_{f_{2n-1}}^{(k)},T_{f_{2n}}^{(k)}]-T_  {\{ f_{2n-1},f_{2n}\} }^{(k)} 
$$
$$
+T_  {\{ f_{2n-1},f_{2n}\} }^{(k)} \bigr ) -
T_{\{ f_1,f_2\} } ^{(k)}...T_ {\{ f_{2n-1},f_{2n}\} }^{(k)}||+O(\frac{1}{k}) . 
$$
This is $O(\frac{1}{k})$. Indeed, within $||.||$  
the term $T_{\{ f_1,f_2\} } ^{(k)}...T_ {\{ f_{2n-1},f_{2n}\} }^{(k)}$ cancels 
and all the other terms are products of factors 
of the form $(ik[T_{f_{2j-1}}^{(k)},T_{f_{2j}}^{(k)}]-T_ {\{ f_{2j-1},f_{2j}\} }^{(k)})$ 
(at least one of these appears) and of the form $T_ {\{ f_{2j-1},f_{2j}\} }^{(k)}$. 
Using the triangle inequality, (\ref{asympcomm}) and Theorem \ref{bmstheorem} (ii), 
we get $O(\frac{1}{k})$. 
Thus, as $k\to\infty$,  
$$
||(ik)^n[T_{f_1}^{(k)},T_{f_2}^{(k)}]...[T_{f_{2n-1}}^{(k)},T_{f_{2n}}^{(k)}]
-T_{\{ f_1,f_2\} ...\{ f_{2n-1},f_{2n}\} }^{(k)}||=O(\frac{1}{k}). 
$$
Exact same proof shows that 
$$
||(ik)^n[T_{f_{\sigma (1)}}^{(k)},T_{f_{\sigma(2) }}^{(k)}]...
[T_{f_{\sigma(2n-1)}}^{(k)},T_{f_{\sigma(2n)}}^{(k)}]
-T_{\{ f_{\sigma(1)},f_{\sigma(2)}\} ...\{ f_{\sigma(2n-1)},f_{\sigma(2n) }\} }^{(k)}||=O(\frac{1}{k}). 
$$
We note that 
$$
T_{\{ f_1,...,f_{2n}\}}^{(k)}=\frac{1}{2^n n!}\sum_{\sigma\in S_{2n}}\sign(\sigma)
T_{\prod_{j=1}^n \{ f_{\sigma (2j-1},f_{\sigma (2j)}\}} ^{(k)} 
$$
(by Lemma \ref{lembrackets}).  
The desired statement now follows from Lemma \ref{lemcomm2} and the triangle inequality. 
$\Box$

The following proposition is similar to Prop. \ref{propcomm}. It   
implies that $\lim_{k\to\infty }||[T_{f_1}^{(k)}, ..., T_{f_{2n}}^{(k)}]||=0$ 
(i.e. $T_{f_1}^{(k)}$, ..., $T_{f_{2n}}^{(k)}$ "Nambu-commute as $k\to\infty$").
\begin{proposition}
\label{commvolform}
For $f_1,...,f_{2n}\in C^{\infty} (M)$
$$
||[T_{f_1}^{(k)}, ..., T_{f_{2n}}^{(k)}]||=O(\frac{1}{k^n})
$$
as $k\to\infty$. 
\end{proposition}
\noindent {\bf Proof.} 
$$
||[T_{f_1}^{(k)}, ..., T_{f_{2n}}^{(k)}]||= 
$$
$$
|| \sideset{}{'} 
\sum _{\sigma\in S_{2n}}\sign(\sigma ) 
[T_{f_{\sigma(1)}}^{(k)}, T_{f_{\sigma(2)}}^{(k)}]
[T_{f_{\sigma(3)}}^{(k)}, T_{f_{\sigma(4)}}^{(k)}]... 
[T_{f_{\sigma(2n-1)}}^{(k)}, T_{f_{\sigma(2n)}}^{(k)}]||\le 
$$
$$
\sideset{}{'} 
\sum  _{\sigma\in S_{2n}}
||[T_{f_{\sigma(1)}}^{(k)}, T_{f_{\sigma(2)}}^{(k)}]||... 
||[T_{f_{\sigma(2n-1)}}^{(k)}, T_{f_{\sigma(2n)}}^{(k)}]||  
$$
which is $O(\frac{1}{k^n})$ by Remark \ref{remarkcomm}. 
$\Box$

\section{Quantization on a hyperk\"ahler manifold}
\label{hyperkahler}

Let $(M,g,J_1,J_2,J_3)$ be a compact connected hyperk\"ahler manifold. Let $4q$ denote the real dimension of $M$.
Denote $\omega_r=g(.,J_r.)$ for $r=1,2,3$.   
The $4$-form
$$
\Omega = \omega_1 \wedge \omega_1 +\omega_2 \wedge \omega_2 + \omega_3 \wedge \omega_3
$$
is $3$-plectic \cite{cantr:99}.
Define the brackets $\{.,.,.,.\} _r$, $\{.,.,.,.\} _{hyp}$ (multilinear maps 
$\bigwedge^4 C^{\infty}(M)\to C^{\infty}(M)$) as follows:
$$
\{f_1,f_2,f_3,f_4\} _r=\{ f_1,f_2\} _r\{ f_3,f_4\} _r- \{ f_1,f_3\} _r\{ f_2,f_4\} _r+  
\{ f_1,f_4\} _r\{ f_2,f_3\} _r ,  
$$
where $\{ .,.\} _r$ is the Poisson bracket on $(M,\omega_r )$, $r=1,2,3$,
$$
\{f_1,f_2,f_3,f_4\} _{hyp} =\sum_{r=1}^{3}\{f_1,f_2,f_3,f_4\} _r .
$$
From the properties of the Poisson bracket it immediately follows that the Leibniz rule is satisfied:
$$
\{f_1,f_2,f_3,f_4f_5\} _{r}= f_4 \{f_1,f_2,f_3,f_5\} _{r}+ \{f_1,f_2,f_3,f_4\} _{r} f_5 
$$
$$
\{f_1,f_2,f_3,f_4f_5\} _{hyp}= f_4 \{f_1,f_2,f_3,f_5\} _{hyp}+ \{f_1,f_2,f_3,f_4\} _{hyp} f_5 .
$$
Therefore $\{ .,.,.,.\} _{r}$, $\{ .,.,.,.\} _{hyp}$  are almost Poisson brackets of order $4$. 

For $q=1$  $\omega_r\wedge \omega_r$ ($r=1,2,3$) and $\Omega$ are volume forms.  
The standard bracket $\{.,.,.,.\} ^{(r)}$ is defined by
$$
df_1\wedge df_2\wedge df_3\wedge df_4=\{f_1,f_2,f_3,f_4\}^{(r)}\frac{1}{2} \omega_r \wedge \omega_r .
$$
From Lemma \ref{lembrackets}, or by a direct calculation (using Darboux theorem, in local coordinates), 
we get: 
\begin{lemma}
For $q=1$ $\{.,.,.,.\} _r$ coincides with $\{.,.,.,.\} ^{(r)}$.  
\end{lemma}
From \cite[Cor. 1 p.106]{gautheron:96}  
it immediately follows that for $q=1$ ($M$ is $4$-dimensional) the Fundamental Identity 
$$
\{f_1,f_2,f_3,\{g_1,g_2,g_3,g_4\}_{hyp}\} _{hyp} =  \{ \{f_1,f_2,f_3,g_1\}_{hyp} ,g_2,g_3,g_4\}_{hyp} +
$$
$$
\{ g_1,\{f_1,f_2,f_3,g_2\}_{hyp}, g_3,g_4\}_{hyp} + 
\{ g_1,g_2,\{f_1,f_2,f_3,g_3\}_{hyp}, g_4\}_{hyp} +
$$
$$
\{ g_1,g_2,g_3,\{f_1,f_2,f_3,g_4\}_{hyp} \}_{hyp}
$$
is satisfied (similarly for $\{.,.,.,.\} _r$). 
For $q>1$ $\{.,.,.,.\} _r$, $\{.,.,.,.\} _{hyp}$ are not necessarily Nambu-Poisson brackets 
(the Fundamental Identity may not be satisfied if $q>1$).  

Assume that the K\"ahler forms $\frac{\omega_1}{2\pi}$, $\frac{\omega_2}{2\pi}$, $\frac{\omega_3}{2\pi}$ are integral.

Let $L_r$ be a holomorphic Hermitian line bundle with curvature of the Chern connection equal to 
$-i\omega_r$, for $r=1,2,3$.
For a positive integer $k$ and $f\in C^\infty (M)$ denote by
$T_{f;r}^{(k)}\in \End(H^0(M,L_r^{\otimes k}))$ the Berezin-Toeplitz operator for $f$.

There are two obvious ways to form a Hilbert space out of three Hilbert spaces 
$H^0(M,L_r^{\otimes k})$ ($r=1,2,3$): by taking direct sum or tensor product. 
Another way to approach this is to say that the vector space of quantization is 
$H^0(M,(L_1\otimes L_2\otimes L_3)^{\otimes k})$, - this would be just the usual Berezin-Toeplitz 
quantization, with the line bundle $L_1\otimes L_2\otimes L_3$. 
Note: in general $H^0(M,(L_1\otimes L_2\otimes L_3)^{\otimes k})$ is not isomorphic 
to $H^0(M,L_1^{\otimes k})\otimes H^0(M,L_2^{\otimes k})\otimes H^0(M,L_3^{\otimes k})$. 

Of course, the hyperk\"ahler structure defines a whole $S^2$ of complex structures (and of K\"ahler forms) on $M$, not just three. A. Uribe pointed out to us that maybe an appropriate notion 
of quantization on a hyperk\"ahler manifold should take into account all $J\in S^2$, 
and should involve an appropriate vector bundle over the twistor space, 
with fibers $H^0(M,L_J^{\otimes k})$. We look forward to seeing his work on this. 

Note that the twistor space of a hyperk\"ahler manifold is not K\"ahler 
(it is generally well-known, see for example \cite{kaledin:99} p. 37, or \cite{huybr:10}), 
so it's not possible to construct a Berezin-Toeplitz quantization on the twistor space. 

\begin{remark}
Denote by $\pi_r : M\times M\times M\to M$ the projection to the $r$-th factor ($r=1,2,3$).  
For sufficiently large $k$ 
$$
H^0(M\times M\times M, (\pi_1^*L_1\otimes \pi_2^*L_2\otimes \pi_3^*L_3)^{\otimes k)})\cong 
\hk
$$
The proof was explained to us by K. Yoshikawa and it goes as follows: 
$$
\dim H^0(M\times M\times M, (\pi_1^*L_1\otimes \pi_2^*L_2\otimes \pi_3^*L_3)^{\otimes k)}=
$$
$$
\int_{M\times M\times M} \td (M\times M\times M)
\ch ((\pi_1^*L_1\otimes \pi_2^*L_2\otimes \pi_3^*L_3)^{\otimes k})=
$$
$$
\int_{M\times M\times M} \pi_1^*\td (M)\pi_2^*\td (M)\pi_3^*\td (M)
\pi_1^*\ch(L_1^{\otimes k}) \pi_2^*\ch (L_2^{\otimes k}) \pi_3^*ch(L_3^{\otimes k})=
$$
$$
\int_M \td (M)\ch(L_1^{\otimes k}) \ \int_M \td (M)\ch(L_2^{\otimes k}) \ 
\int_M \td (M)\ch(L_3^{\otimes k})=
$$
$$
\dim H^0(M,L_1^{\otimes k})\dim H^0(M,L_2^{\otimes k})\dim H^0(M,L_3^{\otimes k})
$$ $\Box$
\end{remark}

In this paper we shall work with functions and structures on $M$, rather than on 
$M\times M\times M$. 

We shall find useful the following statement. 

\begin{proposition}
\label{4functions}
For $f,g,h,t\in C^{\infty} (M)$, $r=1,2,3$, 
$$
||-\frac{k^2}{2}[T_{f;r}^{(k)},T_{g;r}^{(k)},T_{h;r}^{(k)},T_{t;r}^{(k)}]-T_{\{ f,g,h,t\} _r;r}^{(k)}|| = 
O(\frac{1}{k})
$$
as $k\to\infty$.
\end{proposition}
\noindent {\bf Proof.} 
As $k\to\infty$, for $r=1,2,3$, by Theorem \ref{bmstheorem} 
(i) for $f,g\in C^{\infty}(M)$  
\begin{equation}
\label{asympfg}
||ik[T_{f;r}^{(k)},T_{g;r}^{(k)}]-T_{\{ f,g\} _r ;r}^{(k)}||=O(\frac{1}{k}),
\end{equation}
\begin{equation}
\label{asympht}
||ik[T_{h;r}^{(k)},T_{t;r}^{(k)}]-T_{\{ h,t\} _r ;r }^{(k)}||=O(\frac{1}{k}).
\end{equation}
Using Prop. \ref{bmsprop}, we get:
$$
||(ik)^2[T_{f;r}^{(k)},T_{g;r}^{(k)}][T_{h;r}^{(k)},T_{t;r}^{(k)}]-
T_{\{ f,g\} _r \{ h,t\} _r ;r }^{(k)}|| \le 
$$
$$
||(ik)^2[T_{f;r}^{(k)},T_{g;r}^{(k)}][T_{h;r}^{(k)},T_{t;r}^{(k)}]-T_{\{ f,g\} _r ;r }^{(k)}
T_{\{ h,t\} _r ;r }^{(k)}||+
||T_{\{ f,g\} _r ;r }^{(k)}T_{\{ h,t\} _r ;r }^{(k)} -   T_{\{ f,g\} _r \{ h,t\} _r ;r }^{(k)}||=
$$  
$$
||(ik[T_{f;r}^{(k)},T_{g;r}^{(k)}] - T_{\{ f,g\} _r;r }^{(k)} + T_{\{ f,g\} _r;r }^{(k)})
(ik[T_{h;r}^{(k)},T_{t;r}^{(k)}] - T_{\{ h,t\} _r;r }^{(k)} + T_{\{ h,t\} _r;r }^{(k)})
$$
$$
-T_{\{ f,g\} _r;r }^{(k)}T_{\{ h,t\} _r;r }^{(k)}||+O(\frac{1}{k})=
$$
$$
|| (ik[T_{f;r}^{(k)},T_{g;r}^{(k)}] - T_{\{ f,g\} _r;r }^{(k)})(ik[T_{h;r}^{(k)},T_{t;r}^{(k)}] - 
T_{\{ h,t\} _r;r }^{(k)})+
$$
$$
(ik[T_{f;r}^{(k)},T_{g;r}^{(k)}] - T_{\{ f,g\} _r;r }^{(k)})T_{\{ h,t\} _r;r }^{(k)}+
T_{\{ f,g\} _r;r }^{(k)}(ik[T_{h;r}^{(k)},T_{t;r}^{(k)}] - T_{\{ h,t\} _r;r }^{(k)})||+O(\frac{1}{k})
$$
$$
\le || ik([T_{f;r}^{(k)},T_{g;r}^{(k)}] - T_{\{ f,g\} _r;r }^{(k)})|| \ ||ik[T_{h;r}^{(k)},T_{t;r}^{(k)}] - T_{\{ h,t\} _r;r }^{(k)})||+
$$
$$
||ik[T_{f;r}^{(k)},T_{g;r}^{(k)}] - T_{\{ f,g\} _r;r }^{(k)}|| \ ||T_{\{ h,t\} _r;r }^{(k)}||+
||T_{\{ f,g\} _r;r }^{(k)}|| \ ||ik[T_{h;r}^{(k)},T_{t;r}^{(k)}] - T_{\{ h,t\} _r;r }^{(k)}||+O(\frac{1}{k})=
$$
$$
O(\frac{1}{k})O(\frac{1}{k})+|\{ h,t\} _r|_{\infty}O(\frac{1}{k})+|\{ f,g\} _r|_{\infty}O(\frac{1}{k})
+O(\frac{1}{k})=O(\frac{1}{k}).
$$
In the last line we used (\ref{asympfg}), (\ref{asympht}), and applied Theorem \ref{bmstheorem} (ii) twice.  
Similarly we conclude, for $f,h$ and $g,t$:
$$
||(ik)^2[T_{f;r}^{(k)},T_{h;r}^{(k)}][T_{g;r}^{(k)},T_{t;r}^{(k)}]-
T_{\{ f,h\} _r \{ g,t\} _r ;r }^{(k)}||=O(\frac{1}{k}),  
$$
etc. (i.e. we get similar asymptotics for $f,t$ and $g,h$, 
for $h,t$ and $f,g$, for $g,t$ and $f,h$, for $g,h$ and $f,t$). 
Note: 
$$
T_{\{ f,g,h,t\} _r;r }^{(k)}= 
T_{\{ f,g\} _r \{ h,t\} _r ;r }^{(k)}- 
T_{\{ f,h\} _r \{ g,t\} _r ;r }^{(k)}+
T_{\{ f,t\} _r \{ g,h\} _r ;r }^{(k)}. 
$$
Therefore, by (\ref{comm4}) and the triangle inequality, 
$$
||-\frac{k^2}{2}[T_{f;r}^{(k)}, T_{g;r}^{(k)},T_{h;r}^{(k)}, T_{t;r}^{(k)}]-T_{\{ f,g,h,t\} _r;r }^{(k)}||=
O(\frac{1}{k}).
$$ 
$\Box$

\subsection{Direct sum}
\label{directsumgen}

Denote
$$
\hk = H^0(M,L_1^{\otimes k})\oplus H^0(M,L_2^{\otimes k})\oplus H^0(M,L_3^{\otimes k})
$$
(direct sum of Hilbert spaces) 
and
$$
{\bft}_{f}^{(k)}=T_{f;1}^{(k)}\oplus T_{f;2}^{(k)}\oplus T_{f;3}^{(k)},
$$
(${\bft}_{f}^{(k)}$ acts on $\hk$ by ${\bft}_{f}^{(k)}(s_1,s_2,s_3) =(T_{f;1}^{(k)}s_1,T_{f;2}^{(k)}s_2, T_{f;3}^{(k)}s_3)$).
\begin{remark}
Since $||{\bft}_{f}^{(k)}||=\max \{ ||T_{f;1}^{(k)}||,||T_{f;2}^{(k)}||,||T_{f;3}^{(k)}||\}$, 
we immediately have:
\begin{itemize}
\item For $f,g\in C^{\infty}(M)$, as $k\to\infty$,
$$
||ik[\bft_f^{(k)},\bft_g^{(k)}]-\bft_{\{ f,g\} }^{(k)}||=O(\frac{1}{k}), \ 
||[\bft_f^{(k)},\bft_g^{(k)}]||=O(\frac{1}{k})
$$
\item For $f\in C^{\infty}(M)$,  there is a constant $C=C(f)>0$ such that, as $k\to\infty$, 
$$
|f|_{\infty}-\frac{C}{k}\le ||\bft_f^{(k)}||\le |f|_{\infty} . 
$$ 
\item 
For $f_1,...,f_p\in C^{\infty} (M)$ 
$$
||\bft_{f_1}^{(k)}...\bft_{f_p}^{(k)}-\bft_{f_1...f_p}^{(k)}||=O(\frac{1}{k}) 
$$
as $k\to\infty$. 
\end{itemize}
\end{remark}

For $f,g,h,t\in C^{\infty}(M)$ we have:
$$
[\bft_{f}^{(k)},\bft_{g}^{(k)},\bft_{h}^{(k)},\bft_{t}^{(k)}]=\oplus_{r=1}^3 [T_{f;r}^{(k)},T_{g;r}^{(k)},T_{h;r}^{(k)},T_{t;r}^{(k)}] .
$$

\begin{theorem}
\label{thhyperk}
For $f,g,h,t\in C^{\infty} (M)$
$$
||-\frac{k^2}{2}[\bft_{f}^{(k)},\bft_{g}^{(k)},\bft_{h}^{(k)},\bft_{t}^{(k)}]-\oplus_{r=1}^3T_{\{ f,g,h,t\} _r;r}^{(k)}|| = 
O(\frac{1}{k})
$$
as $k\to\infty$.
\end{theorem}
\noindent {\bf Proof.} 
Using Proposition \ref{4functions}, we get:
$$
||-\frac{k^2}{2}[\bft_f^{(k)}, \bft_g^{(k)},\bft_h^{(k)}, \bft_t^{(k)}]-\oplus_{r=1}^3 T_{\{ f,g,h,t\} _r;r }^{(k)}||=
$$
$$
\max _{1\le r\le 3} 
||-\frac{k^2}{2}[T_{f;r}^{(k)}, T_{g;r}^{(k)},T_{h;r}^{(k)}, T_{t;r}^{(k)}]-T_{\{ f,g,h,t\} _r;r }^{(k)}||=
O(\frac{1}{k}).
$$
$\Box$   

The following proposition is similar to Prop. \ref{propcomm}. 
It implies that $\bft_f^{(k)}$, $\bft_g^{(k)}$, $\bft_h^{(k)}$, $\bft_t^{(k)}$ "Nambu-commute as $k\to\infty$".
\begin{proposition}
\label{commdim4}
For $f_1,f_2,f_3,f_4\in C^{\infty} (M)$
$$
||[\bft _{f_1}^{(k)}, \bft_{f_2}^{(k)}, \bft_{f_3}^{(k)}, \bft_{f_4}^{(k)}]||=O(\frac{1}{k^2}) 
$$
as $k\to\infty$. 
\end{proposition}
\noindent {\bf Proof.} 
$$
||[\bft _{f_1}^{(k)}, \bft _{f_2}^{(k)}, \bft _{f_3}^{(k)}, \bft _{f_4}^{(k)}]||=
\max _{1\le r\le 3}||[T_{f_1;r}^{(k)}, T_{f_2;r}^{(k)}, T_{f_3;r}^{(k)}, T_{f_4;r}^{(k)}]||= 
$$
$$
\max _{1\le r\le 3}||\sideset{}{'} \sum  _{\sigma\in S_4}\sign(\sigma) 
[T_{f_{\sigma(1)};r}^{(k)}, T_{f_{\sigma(2)};r}^{(k)}]
[T_{f_{\sigma(3)};r}^{(k)}, T_{f_{\sigma(4)};r}^{(k)}]||\le 
$$
$$
\max _{1\le r\le 3}\sideset{}{'} \sum_{\sigma\in S_4}
||[T_{f_{\sigma(1)};r}^{(k)}, T_{f_{\sigma(2)};r}^{(k)}]|| \ 
||[T_{f_{\sigma(3)};r}^{(k)}, T_{f_{\sigma(4)};r}^{(k)}]||. 
$$
By Remark \ref{remarkcomm} it is $O(\frac{1}{k^2})$. 
$\Box$ 
 
\subsection{Direct sum: dimension 4}
\label{directsum4}

To discuss the correspondence between the the bracket on functions 
and the generalized commutator (as $k\to\infty$) in the hyperk\"ahler case: 
we showed (Theorem \ref{thhyperk}) 
that for a hyperk\"ahler manifold $M$ 
of arbitrary dimension and smooth functions $f,g,h,t$ on $M$ 
$[\bft _{f}^{(k)},\bft _{g}^{(k)},\bft _{h}^{(k)},\bft _{t}^{(k)}]$ 
is asymptotic to 
$$
\begin{pmatrix} T_{\{ f,g,h,t\} _1;1}^{(k)} & & \\ 
& T_{\{ f,g,h,t\} _2;2}^{(k)} & \\
& & T_{\{ f,g,h,t\} _3;3}^{(k)}
 \end{pmatrix} ,
$$
{\bf not} to 
$$
\bft _{\{ f,g,h,t\}_{hyp}}^{(k)}=
\begin{pmatrix} T_{\{ f,g,h,t\}_{hyp} ;1}^{(k)} & & \\ 
& T_{\{ f,g,h,t\}_{hyp} ;2}^{(k)} & \\
& & T_{\{ f,g,h,t\}_{hyp} ;3}^{(k)}
 \end{pmatrix} .
$$ 
To clarify, we have obtained an asymptotic relation between a map 
$$
\sideset{}{^4} 
\bigwedge C^{\infty}(M) \to C^{\infty}(M)\times C^{\infty}(M) \times C^{\infty}(M) 
$$
$$
f,g,h,t \mapsto (\{ f,g,h,t\} _1, \{ f,g,h,t\} _2, \{ f,g,h,t\} _3)
$$ 
and the Nambu generalized commutator $[.,.,.,.]$. 
It is not the same as a correspondence between 
$\{ .,.,.,.\}_{hyp} : \bigwedge ^4 C^{\infty}(M) \to C^{\infty}(M)$ and $[.,.,.,.]$. 

From now on $M$ will be of real dimension $4$ (hence $M$ is isomorphic 
to a K3-surface or a torus \cite{besse:87} 14.22). In this case we get 
Theorem \ref{thdim4} below, and in the case 
when $M$ is a $4$-torus with three standard linear complex structures 
(Example \ref{R4} below) - 
we get that $[\bft _{f}^{(k)},\bft _{g}^{(k)},\bft _{h}^{(k)},\bft _{t}^{(k)}]$ 
is asymptotic to $\bft _{\{ f,g,h,t\}_{hyp}}^{(k)}$. 

We have: for $r=1,2,3$  
$$
\Omega = \frac{\mu _r}{2}\omega_r\wedge\omega_r, 
$$
where $\mu_r$ is a smooth non-vanishing function on $M$. 
Denote by $\{ .,.,.,.\}$ the Nambu-Poisson bracket defined by 
$$
df_1\wedge df_2\wedge df_3\wedge df_4=\{ f_1,f_2,f_3,f_4\} \Omega .
$$  
Therefore 
$$
\{ f_1,f_2,f_3,f_4\} _r=\{ f_1,f_2,f_3,f_4\} ^{(r)}=\mu_r \{ f_1,f_2,f_3,f_4\} . 
$$
Denote 
$$
\bft _{\mu}^{(k)}=\begin{pmatrix}
T_{\mu_1;1}^{(k)} & & \\
& T_{\mu_2;2}^{(k)} & \\
& & T_{\mu_3;3}^{(k)}
\end{pmatrix} .
$$
The following theorem shows that 
$[\bft _{f_1}^{(k)},\bft _{f_2}^{(k)},\bft _{f_3}^{(k)},\bft _{f_4}^{(k)}]$ 
is asymptotic to $\bft _{\{ f_1,f_2,f_3,f_4\} }^{(k)}\bft _{\mu}^{(k)}$. 
\begin{theorem}
\label{thdim4}
For $f,g,h,t\in C^{\infty}(M)$ 
$$
||-\frac{k^2}{2} [ \bft _f^{(k)},\bft _g^{(k)},\bft _h^{(k)},\bft _t^{(k)}]
-\bft _{\{ f,g,h,t\} }^{(k)}\bft _{\mu}^{(k)}||
=O(\frac{1}{k})
$$
as $k\to\infty$. 
\end{theorem} 
\noindent {\bf Proof.} 
For $r=1,2,3$ the same argument as in the proof of Proposition \ref{4functions} gives: 
\begin{equation}
\label{asympt}
||-\frac{k^2}{2} [ T_{f;r}^{(k)},T_{g;r}^{(k)},T_{h;r}^{(k)},T_{t;r}^{(k)}]-
T_{\{ f,g,h,t\} _r ;r }^{(k)}|| 
=O(\frac{1}{k})
\end{equation}
We have: 
$$
||-\frac{k^2}{2} [ T_{f;r}^{(k)},T_{g;r}^{(k)},T_{h;r}^{(k)},T_{t;r}^{(k)}]-
T_{\{ f,g,h,t\} ;r }^{(k)}T_{\mu _r;r}^{(k)}|| \le 
$$
$$
||-\frac{k^2}{2} [ T_{f;r}^{(k)},T_{g;r}^{(k)},T_{h;r}^{(k)},T_{t;r}^{(k)}]-
T_{\{ f,g,h,t\} \mu _r;r}^{(k)}|| +
||T_{\{ f,g,h,t\} \mu _r;r}^{(k)} - T_{\{ f,g,h,t\} ;r }^{(k)}T_{\mu _r;r}^{(k)}||. 
$$
This is $O(\frac{1}{k})$ by (\ref{asympt}) and Prop. \ref{bmsprop}. 
Hence 
$$
||-\frac{k^2}{2} [ \bft _f^{(k)},\bft _g^{(k)},\bft _h^{(k)},\bft _t^{(k)}]
-\bft _{\{ f,g,h,t\} }^{(k)}\bft _{\mu}^{(k)}||
$$
$$
=\max_{1\le r\le 3} 
||-\frac{k^2}{2} [ T_{f;r}^{(k)},T_{g;r}^{(k)},T_{h;r}^{(k)},T_{t;r}^{(k)}]-
T_{\{ f,g,h,t\}  ;r }^{(k)}T_{\mu _r ;r}^{(k)}||=O(\frac{1}{k}).
$$
$\Box$
\begin{example} 
\label{R4}
Denote $\tilde{M}=\R^4$, with coordinates $x_1$, $x_2$, $x_3$, $x_4$, 
and equipped with three (linear) complex structures 
$$
J_1=\begin{pmatrix}
0 & -1 & & \\
1 & 0 & & \\
& & 0 & -1 \\
& & 1 & 0 
\end{pmatrix}, \  
J_2=\begin{pmatrix}
 &  & -1 & 0 \\
 &  & 0 &  1 \\
1 & 0 &  &  \\
0 & -1 &  &  
\end{pmatrix}, \ 
J_3=\begin{pmatrix}
 & & & -1 \\
 &  & -1 & \\
& 1 &  &  \\
1 & & &  
\end{pmatrix} .
$$
We have: $J_1J_2=J_3$ and, of course, $J_1^2=J_2^2=J_3^2=-I$. 

Note: if we regard $\tilde{M}$ as the one-dimensional quaternionic vector space, 
with basis 1, {\bf{i, j, k}} ({\bf i}$^2$={\bf j}$^2$={\bf k}$^2$=$-1$, {\bf ij=k}), then 
$J_1$, $J_2$, $J_3$ correspond to left multiplication by {\bf{i, j, k}} respectively. 

For the standard Riemannian metric on $\tilde{M}$, with the metric tensor $g=I$, 
the symplectic forms are as follows: 
$$
\omega_1 = dx_1\wedge dx_2 +dx_3\wedge dx_4,
$$
$$
\omega_2 = dx_1\wedge dx_3 -dx_2\wedge dx_4,
$$
$$
\omega_3 = dx_1\wedge dx_4 +dx_2\wedge dx_3. 
$$
For $r=1,2,3$ 
$$
\frac{1}{2}\omega_r\wedge \omega_r=dx_1\wedge dx_2\wedge dx_3\wedge dx_4,
$$
$$
\Omega =\sum_{r=1}^3 \omega_r\wedge \omega_r = 6 dx_1\wedge dx_2\wedge dx_3\wedge dx_4. 
$$
Everything is $\Z^4$-invariant and $g$, $J_1$, $J_2$, $J_3$, 
$\omega_1$, $\omega_2$, $\omega_3$, $\Omega$  descend to $M=\tilde{M}/\Z^4$. 
We get: $\mu_1=\mu_2=\mu_3=6$ and 
$$
6\{ .,.,.,.\} = \{ .,.,.,.\} _r=\{ .,.,.,.\} ^{(r)}=\frac{1}{3}\{ .,.,.,.\} _{hyp}. 
$$
Theorem \ref{thdim4} gives:  
for $f,g,h,t\in C^{\infty}(M)$ 
\begin{equation}
\label{asymptorus}
||-c k^2 [ \bft _f^{(k)},\bft _g^{(k)},\bft _h^{(k)},\bft _t^{(k)}]
-\bft _{\{ f,g,h,t\} _{hyp}}^{(k)}||
=O(\frac{1}{k})
\end{equation}
as $k\to\infty$, where $c$ is a positive constant. 
\end{example}

\subsection{Tensor product}
\label{tensorproduct}

Denote
$$
\hk = H^0(M,L_1^{\otimes k})\otimes H^0(M,L_2^{\otimes k})\otimes H^0(M,L_3^{\otimes k})
$$
(tensor product of Hilbert spaces) 
and
$$
\bbt ^{(k)}_f=T_{f;1}^{(k)}\otimes T_{f;2}^{(k)}\otimes T_{f;3}^{(k)},
$$
($\bbt _{f}^{(k)}(s_1\otimes s_2\otimes s_3) =T_{f;1}^{(k)}s_1\otimes T_{f;2}^{(k)}s_2\otimes 
T_{f;3}^{(k)}s_3$ and the action extends to $\hk$ by linearity, 
also note: $||\bbt_{f}^{(k)}||=||T_{f;1}^{(k)}|| \ ||T_{f;2}^{(k)}|| \ ||T_{f;3}^{(k)}||$).

In the proofs below we shall need the following elementary statement. 
\begin{lemma}
\label{tensorlem} 
If $M_j$, $N_j$ are linear operators on a finite dimensional Hilbert space $V_j$ 
($j=1,2,3$), then 
$$
||M_1\otimes M_2\otimes M_3-N_1\otimes N_2\otimes N_3||\le 
||M_1-N_1|| \ ||M_2-N_2|| \  ||M_3-N_3||+
$$
$$
||M_1-N_1|| \ || M_2|| \  ||N_3||+
||M_1|| \  ||N_2|| \ ||M_3-N_3||+||N_1|| \ ||M_2-N_2|| \  ||M_3||
$$
\end{lemma}
\noindent {\bf Proof.} This immediately follows from the equality 
$$
(M_1-N_1)\otimes (M_2-N_2)\otimes (M_3-N_3)=M_1\otimes M_2\otimes M_3-N_1\otimes N_2\otimes N_3-
$$
$$
(M_1-N_1)\otimes M_2\otimes N_3-M_1\otimes N_2\otimes (M_3-N_3)-N_1\otimes (M_2-N_2)\otimes M_3
$$
$\Box$

We also note the following identity for tensor products of operators: 
\begin{equation}
\label{tensorid3}
[A_1\otimes A_2 \otimes A_3, B_1\otimes B_2 \otimes B_3]=
[A_1,B_1]\otimes [A_2,B_2] \otimes [A_3,B_3] + 
\end{equation}
$$
[A_1,B_1]\otimes B_2A_2\otimes A_3B_3+ 
A_1B_1 \otimes [A_2,B_2]\otimes B_3A_3 + B_1A_1 \otimes A_2B_2 \otimes [A_3,B_3]. 
$$
\begin{remark}  \ 

\begin{itemize}
\item For $f\in C^{\infty}(M)$,  there is a constant $C=C(f)>0$ such that, as $k\to\infty$,
$$
(|f|_{\infty}-\frac{C}{k})^3\le ||\bbt_f^{(k)}||\le (|f|_{\infty})^3 . 
$$ 
\item 
For $f_1,...,f_p\in C^{\infty} (M)$ 
$$
||\bbt_{f_1}^{(k)}...\bbt_{f_p}^{(k)}-\bbt_{f_1...f_p}^{(k)}||=O(\frac{1}{k}) 
$$
as $k\to\infty$. 
\end{itemize}
The last statement holds for $p=2$ by Lemma \ref{tensorlem}, Theorem \ref{bmstheorem} 
and Prop. \ref{bmsprop}. It follows for arbitrary $p$ by induction.  
\end{remark}

\begin{proposition}
\label{3comms}
For $f,g\in C^{\infty} (M)$ 
$$
||(ik)^3[T_{f;1}^{(k)},T_{g;1}^{(k)}]\otimes [T_{f;2}^{(k)},T_{g;2}^{(k)}]
\otimes [T_{f;3}^{(k)},T_{g;3}^{(k)}]-
T_{\{ f,g\} _1;1}^{(k)}\otimes T_{\{ f,g\} _2;2}^{(k)}\otimes 
T_{\{ f,g\} _3;3}^{(k)}||=O(\frac{1}{k})
$$
as $k\to\infty$. 
\end{proposition}
\noindent {\bf Proof.} This follows from Lemma \ref{tensorlem}, Theorem \ref{bmstheorem} 
and Remark \ref{remarkcomm}.  $\Box$

\begin{proposition}
\label{tensorprop}
For $f,g\in C^{\infty} (M)$ 
$$
||ik[\bbt_f^{(k)},\bbt_g^{(k)}]-
(T_{\{ f,g\} _1;1}^{(k)}\otimes T_{fg;2}^{(k)}\otimes 
T_{fg;3}^{(k)}+
T_{fg;1}^{(k)}\otimes T_{\{ f,g\} _2;2}^{(k)}\otimes 
T_{fg;3}^{(k)}+
$$
$$
T_{fg;1}^{(k)}\otimes T_{fg;2}^{(k)}\otimes 
T_{\{ f,g\} _3;3}^{(k)})
||=O(\frac{1}{k})
$$
as $k\to\infty$. 
\end{proposition}

\noindent {\bf Proof.} Using (\ref{tensorid3}), we get: 
$$
||ik[\bbt_f^{(k)},\bbt_g^{(k)}]-
(T_{\{ f,g\} _1;1}^{(k)}\otimes T_{fg;2}^{(k)}\otimes 
T_{fg;3}^{(k)}+
T_{fg;1}^{(k)}\otimes T_{\{ f,g\} _2;2}^{(k)}\otimes 
T_{fg;3}^{(k)}+
$$
$$
T_{fg;1}^{(k)}\otimes T_{fg;2}^{(k)}\otimes 
T_{\{ f,g\} _3;3}^{(k)})||\le 
$$
$$
||ik[T_{f;1}^{(k)},T_{g;1}^{(k)}]\otimes T_{g;2}^{(k)}T_{f;2}^{(k)}\otimes 
T_{f;3}^{(k)}T_{g;3}^{(k)}-
T_{\{ f,g\} _1;1}^{(k)}\otimes T_{fg;2}^{(k)}\otimes 
T_{fg;3}^{(k)}||+
$$
$$
||T_{f;1}^{(k)}T_{g;1}^{(k)}\otimes ik[T_{f;2}^{(k)},T_{g;2}^{(k)}]\otimes 
T_{g;3}^{(k)}T_{f;3}^{(k)}-
T_{fg;1}^{(k)}\otimes T_{\{ f,g\} _2;2}^{(k)}\otimes 
T_{fg;3}^{(k)}||+
$$
$$
||T_{g;1}^{(k)}T_{f;1}^{(k)}\otimes T_{f;2}^{(k)}T_{g;2}^{(k)}\otimes 
ik[T_{f;3}^{(k)},T_{g;3}^{(k)}]-
T_{fg;1}^{(k)}\otimes T_{fg;2}^{(k)}\otimes 
T_{\{ f,g\} _3;3}^{(k)}||+
$$
$$
k||[T_{f;1}^{(k)},T_{g;1}^{(k)}]\otimes [T_{f;2}^{(k)},T_{g;2}^{(k)}]\otimes 
[T_{f;3}^{(k)},T_{g;3}^{(k)}]||. 
$$
Each of the first three terms is $O(\frac{1}{k})$ by Lemma \ref{tensorlem},  
Theorem \ref{bmstheorem}, Prop. \ref{bmsprop} and Remark \ref{remarkcomm}. 
The last term is $O(\frac{1}{k^2})$ by Remark \ref{remarkcomm}. 
$\Box$

\begin{corollary}
\label{tensorcorcomm}
For $f,g\in C^{\infty}(M)$
$$
||[\bbt _f^{(k)},\bbt _g^{(k)}]||=O(\frac{1}{k}) 
$$
as $k\to\infty$. 
\end{corollary}
\noindent {\bf Proof.} Follows from Proposition \ref{tensorprop} and Theorem \ref{bmstheorem}(ii)  
by triangle inequality. $\Box$
\begin{corollary}
\label{tensorcorcomm2}
For $f,g,h,t\in C^{\infty}(M)$
$$
||[\bbt _f^{(k)},\bbt _g^{(k)},\bbt _h^{(k)},\bbt _t^{(k)}]||=O(\frac{1}{k^2}) 
$$
as $k\to\infty$. 
\end{corollary}
\noindent {\bf Proof.} Follows from equality (\ref{comm4}) and Corollary \ref{tensorcorcomm} 
by triangle inequality. $\Box$

\begin{proposition}
\label{tensorprop4}
For $f,g,h,t\in C^{\infty}(M)$ 
$$
||-\frac{k^6}{8}[T _{f;1}^{(k)},T_{g;1}^{(k)},T_{h;1}^{(k)},T_{t;1}^{(k)}]\otimes 
[T _{f;2}^{(k)},T_{g;2}^{(k)},T_{h;2}^{(k)},T_{t;2}^{(k)}]\otimes 
[T _{f;3}^{(k)},T_{g;3}^{(k)},T_{h;3}^{(k)},T_{t;3}^{(k)}]- 
$$
$$
T_{\{ f,g,h,t\} _1 ;1}^{(k)}\otimes 
T_{\{ f,g,h,t\} _2 ;2}^{(k)}\otimes 
T_{\{ f,g,h,t\} _3 ;3}^{(k)}
||=O(\frac{1}{k})
$$
as $k\to\infty$.
\end{proposition}
\noindent {\bf Proof.} 
For $r=1,2,3$ 
$$
||[T _{f;r}^{(k)},T_{g;r}^{(k)},T_{h;r}^{(k)},T_{t;r}^{(k)}]||=O(\frac{1}{k^2})
$$
as $k\to\infty$ (this follows by triangle inequality from (\ref{comm4}) and Remark \ref{remarkcomm}). 
The statement now 
follows from Lemma \ref{tensorlem}, Proposition \ref{4functions} and  
Theorem \ref{bmstheorem} (ii).  
$\Box$

It is natural to ask about asymptotics of 
$[\bbt _f^{(k)},\bbt _g^{(k)},\bbt _h^{(k)},\bbt _t^{(k)}]$ for given 
$f,g,h,t \in C^{\infty} (M)$. 
Proposition \ref{tensorprop} dictates the following very technical statement. 

\begin{theorem}
\label{tensorth}
For $f_1,f_2,f_3,f_4\in C^{\infty}(M)$ 
$$
||-\frac{k^2}{2}[\bbt _{f_1}^{(k)},\bbt _{f_2}^{(k)},\bbt _{f_3}^{(k)},\bbt _{f_4}^{(k)}]- 
{\mathbb{W}}_{f_1,f_2,f_3,f_4 }^{(k)}||=O(\frac{1}{k})
$$
as $k\to\infty$, where 
$$
{\mathbb{W}}_{f_1,f_2,f_3,f_4 }^{(k)}=
T_{ \{ f_1,f_2,f_3,f_4\} _1 ;1}^{(k)}\otimes T_{ f_1f_2f_3f_4;2}^{(k)}
\otimes T_{ f_1f_2f_3f_4 ;3}^{(k)}+
$$
$$
T_{ f_1f_2f_3f_4 ;1}^{(k)}\otimes T_{ \{ f_1,f_2,f_3,f_4\} _2 ;2}^{(k)}
\otimes T_{ f_1f_2f_3f_4 ;3}^{(k)}+
T_{f_1f_2f_3f_4 ;1}^{(k)}\otimes T_{f_1f_2f_3f_4 ;2}^{(k)}
\otimes T_{ \{ f_1,f_2,f_3,f_4\} _3 ;3}^{(k)}+
$$
$$
\sum_{\substack{{(i,j,m,l)=(1,2,3,4),}\\ 
{(1,3,2,4),(1,4,2,3)}} } \sign (i,j,m,l) \Bigl [ T_{f_if_j\{ f_m,f_l\} _1;1}^{(k)} \otimes 
( T_{f_mf_l\{ f_i,f_j\} _2;2}^{(k)} \otimes T_{f_if_jf_mf_l;3}^{(k)} +
$$
$$
T_{f_if_jf_mf_l;2}^{(k)} \otimes   T_{f_mf_l\{ f_i,f_j\} _3;3}^{(k)} )+
T_{f_mf_l\{ f_i,f_j\} _1;1}^{(k)}\otimes 
$$ 
$$
(T_{f_if_j\{ f_m,f_l\} _2;2}^{(k)} \otimes T_{f_if_jf_mf_l;3}^{(k)}+
T_{f_if_jf_mf_l;2}^{(k)} \otimes  T_{f_if_j\{ f_m,f_l\} _3;3}^{(k)})+
$$ 
$$
T_{f_if_jf_mf_l;1}^{(k)} \otimes (
T_{f_if_j\{ f_m,f_l\} _2;2}^{(k)} \otimes T_{f_mf_l\{ f_i,f_j\} _3;3}^{(k)}+
T_{f_mf_l\{ f_i,f_j\} _2;2}^{(k)}\otimes T_{f_if_j\{ f_m,f_l\} _3;3}^{(k)}) \Bigr ] .
$$

\end{theorem}
\noindent {\bf Proof.} 
First, we observe: as $k\to\infty$ 
$$
||(ik)^2[\bbt _{f_i}^{(k)},\bbt _{f_j}^{(k)}][\bbt _{f_m}^{(k)},\bbt _{f_l}^{(k)}]-
(T_{\{ f_i,f_j\} _1;1}^{(k)}\otimes T_{f_i f_j;2}^{(k)}\otimes 
T_{f_i f_j;3}^{(k)}+
$$
$$
T_{f_i f_j;1}^{(k)}\otimes T_{\{ f_i,f_j\} _2;2}^{(k)}\otimes 
T_{f_i f_j;3}^{(k)}+
T_{f_i f_j;1}^{(k)}\otimes T_{f_i f_j;2}^{(k)}\otimes 
T_{\{ f_i,f_j\} _3;3}^{(k)})
$$
$$
(T_{\{ f_m,f_l\} _1;1}^{(k)}\otimes T_{f_m f_l;2}^{(k)}\otimes 
T_{f_m f_l;3}^{(k)}+
T_{f_m f_l;1}^{(k)}\otimes T_{\{ f_m,f_l\} _2;2}^{(k)}\otimes 
T_{f_m f_l;3}^{(k)}+
$$
$$
T_{f_m f_l;1}^{(k)}\otimes T_{f_m f_l;2}^{(k)}\otimes 
T_{\{ f_m,f_l\} _3;3}^{(k)})||=O(\frac{1}{k}).
$$
This follows from the elementary inequality 
$$
||M_1M_2-N_1N_2||=||M_1M_2-M_2N_1+M_2N_1-N_1N_2||\le 
$$
$$
||M_2||||M_1-N_1||+||N_1||||M_2-N_2||
$$
by setting 
$$
M_1=ik[\bbt _{f_i}^{(k)},\bbt _{f_j}^{(k)}], \ M_2=ik[\bbt _{f_m}^{(k)},\bbt _{f_l}^{(k)}]
$$
$$
N_1=T_{\{ f_i,f_j\} _1;1}^{(k)}\otimes T_{f_i f_j;2}^{(k)}\otimes 
T_{f_i f_j;3}^{(k)}+
T_{f_i f_j;1}^{(k)}\otimes T_{\{ f_i,f_j\} _2;2}^{(k)}\otimes 
T_{f_i f_j;3}^{(k)}+
$$
$$
T_{f_i f_j;1}^{(k)}\otimes T_{f_i f_j;2}^{(k)}\otimes 
T_{\{ f_i,f_j\} _3;3}^{(k)}, 
$$
$$
N_2=T_{\{ f_m,f_l\} _1;1}^{(k)}\otimes T_{f_m f_l;2}^{(k)}\otimes 
T_{f_m f_l;3}^{(k)}+
T_{f_m f_l;1}^{(k)}\otimes T_{\{ f_m,f_l\} _2;2}^{(k)}\otimes 
T_{f_m f_l;3}^{(k)}+
$$
$$
T_{f_m f_l;1}^{(k)}\otimes T_{f_m f_l;2}^{(k)}\otimes 
T_{\{ f_m,f_l\} _3;3}^{(k)},   
$$
with the use of Theorem \ref{bmstheorem}(ii), Prop. \ref{tensorprop} and Cor. \ref{tensorcorcomm}.   
Next, using Lemma \ref{tensorlem}, Theorem \ref{bmstheorem}(ii) 
and Prop. \ref{bmsprop}, we get: 
$$
||-k^2[\bbt _{f_i}^{(k)},\bbt _{f_j}^{(k)}][\bbt _{f_m}^{(k)},\bbt _{f_l}^{(k)}]-\Bigl [ 
T_{ \{ f_i,f_j\}_1 \{ f_m,f_l\} _1 ;1}^{(k)}\otimes T_{ f_if_jf_mf_l;2}^{(k)}
\otimes T_{ f_if_jf_mf_l ;3}^{(k)}+
$$
$$
T_{ f_if_jf_mf_l ;1}^{(k)}\otimes T_{ \{ f_i,f_j\} _2 \{ f_m,f_l\} _2 ;2}^{(k)}
\otimes T_{ f_if_jf_mf_l ;3}^{(k)}+
T_{f_if_jf_mf_l ;1}^{(k)}\otimes T_{f_if_jf_mf_l ;2}^{(k)}
\otimes T_{ \{ f_i,f_j\} _3 \{ f_m,f_l\} _3 ;3}^{(k)}+
$$
$$
T_{f_if_j\{ f_m,f_l\} _1;1}^{(k)} \otimes 
( T_{f_mf_l\{ f_i,f_j\} _2;2}^{(k)} \otimes T_{f_if_jf_mf_l;3}^{(k)} +
T_{f_if_jf_mf_l;2}^{(k)} \otimes   T_{f_mf_l\{ f_i,f_j\} _3;3}^{(k)} )+
$$
$$
T_{f_mf_l\{ f_i,f_j\} _1;1}^{(k)}\otimes 
(T_{f_if_j\{ f_m,f_l\} _2;2}^{(k)} \otimes T_{f_if_jf_mf_l;3}^{(k)}+
T_{f_if_jf_mf_l;2}^{(k)} \otimes  T_{f_if_j\{ f_m,f_l\} _3;3}^{(k)})+
$$ 
$$
T_{f_if_jf_mf_l;1}^{(k)} \otimes (
T_{f_if_j\{ f_m,f_l\} _2;2}^{(k)} \otimes T_{f_mf_l\{ f_i,f_j\} _3;3}^{(k)}+
T_{f_mf_l\{ f_i,f_j\} _2;2}^{(k)}\otimes T_{f_if_j\{ f_m,f_l\} _3;3}^{(k)}) \Bigr ] ||=O(\frac{1}{k}).
$$

After that we note: 
$$
[\bbt _{f_1}^{(k)},\bbt _{f_2}^{(k)},\bbt _{f_3}^{(k)},\bbt _{f_4}^{(k)}]=
\sum_{\substack{{(i,j,m,l)=}\\ 
{(1,2,3,4),(1,3,2,4),(1,4,2,3)}\\
{(3,4,1,2),(2,4,1,3),(2,3,1,4)}
} 
} \sign (i,j,m,l)[\bbt _{f_i}^{(k)},\bbt _{f_j}^{(k)}][\bbt _{f_m}^{(k)},\bbt _{f_l}^{(k)}]
$$
(see (\ref{comm4})). Taking the sum, we get: 
$$
||-k^2[\bbt _{f_1}^{(k)},\bbt _{f_2}^{(k)},\bbt _{f_3}^{(k)},\bbt _{f_4}^{(k)}]- 
2{\mathbb{W}}_{f_1,f_2,f_3,f_4 }^{(k)}||=O(\frac{1}{k}).
$$
$\Box$


\begin{thebibliography}{999}


\bibitem[A1]{andersen:06}J. Andersen. 
\newblock {\em Asymptotic faithfulness of the quantum SU(n) representations of the mapping class groups.}  
\newblock Ann. Math. (2) 163 (2006), no. 1, 347-368. 

\bibitem[A2]{andersen:10}J. Andersen. 
\newblock {\em Toeplitz operators and Hitchin's projectively flat connection.}  
\newblock In {\it The many facets of geometry}, 177-209, Oxford Univ. Press, Oxford, 2010. 


\bibitem[A]{artin:91}M. Artin.  
\newblock {\em Algebra.} 
\newblock Prentice Hall, 1991. 

\bibitem[AI]{azcar:10}J. A. de Azc\'arraga, J. M. Izquierdo.  
\newblock {\em n-ary algebras: a review with applications.}
\newblock 2010 J. Phys. A: Math. Theor. 43, 293001. 

\bibitem[APP]{azcar:96}J. A. de Azc\'arraga,  A. M. Perelomov, J.C. P\'erez Bueno. 
\newblock {\em New generalized Poisson structures.} 
\newblock J. Phys. A 29 (1996), no. 7, L151-L157. 

\bibitem[BHR]{baez:10}J. Baez, A. Hoffnung, C. Rogers. 
\newblock {\em Categorified symplectic geometry and the classical string.} 
\newblock Comm. Math. Phys. 293 (2010), no. 3, 701-725. 

\bibitem[BR]{baezrog:10}J. Baez, C. Rogers.
\newblock {\em Categorified symplectic geometry and the string Lie 2-algebra.} 
\newblock Homology, Homotopy Appl. 12 (2010), no. 1, 221-236. 

\bibitem[BMMP]{barron:14}T. Barron, X. Ma, G. Marinescu, M. Pinsonnault. 
\newblock {\em Semi-classical properties of Berezin-Toeplitz operators with $C^k$-symbol.}
\newblock J. Math. Phys. 55, 042108 (2014)

\bibitem[BF]{bayen:75}F. Bayen, M. Flato. 
\newblock {\em Remarks concerning Nambu's generalized mechanics.} 
\newblock Phys. Rev. D (3) 11 (1975), 3049-3053. 

\bibitem[Ber]{berezin:74}F. Berezin.  
\newblock {\em Quantization.} 
\newblock Math. USSR-Izv. 8 (1974), no. 5, 1109-1165 (1975). 

\bibitem[Bes]{besse:87}A. Besse. 
\newblock {\em Einstein manifolds.}
\newblock Springer-Verlag, Berlin, 1987. 

\bibitem[BMS]{bordemann:94}M. Bordemann, E. Meinrenken, M. Schlichenmaier.  
\newblock {\em Toeplitz quantization of K\"ahler manifolds and 
gl(N), $N\to\infty$ limits. }
\newblock Comm. Math. Phys. 165 (1994), no. 2, 281-296. 

\bibitem[BG]{boutet:81}L. Boutet de Monvel, V. Guillemin.  
\newblock {\em The spectral theory of Toeplitz operators.}
\newblock Annals of Mathematics Studies, 99. 
\newblock Princeton University Press, Princeton, NJ; 
\newblock University of Tokyo Press, Tokyo, 1981.

\bibitem[Br]{bremner:98}M. Bremner. 
\newblock {\em Identities for the ternary commutator.} 
\newblock J. Algebra 206 (1998), no. 2, 615-623. 

\bibitem[BCI]{bur:13}H. Bursztyn, A. Cabrera, D. Iglesias. 
\newblock {\em Multisymplectic geometry and Lie groupoids.} 
\newblock arXiv:1312.6436 [math.SG], to appear in 
\newblock {\it Geometry, Mechanics and Dynamics: The Legacy of Jerry Marsden}, 
\newblock Fields Institute Communications Series

\bibitem[CIL]{cantr:99}F. Cantrijn, A. Ibort, M. de Le\'on. 
\newblock {\em On the geometry of multisymplectic manifolds.} 
\newblock J. Austral. Math. Soc. Ser. A 66 (1999), no. 3, 303-330. 

\bibitem[CT]{chat:96}R. Chatterjee, L. Takhtajan. 
\newblock {\em Aspects of classical and quantum Nambu mechanics.}
\newblock Lett. Math. Phys. 37 (1996), no. 4, 475-482. 


\bibitem[CJM]{curtr:09}T. Curtright, X. Jin, L. Mezincescu. 
\newblock {\em Multi-operator brackets acting thrice.} 
\newblock J. Phys. A 42 (2009), no. 46, 462001, 6 pp. 

\bibitem[CZ1]{curtr:03}T. Curtright, C. Zachos. 
\newblock {\em Classical and quantum Nambu mechanics.} 
\newblock Phys. Rev. D. 68 (2003) 085001, 1-29. 

\bibitem[CZ2]{curtr:04}T. Curtright, C. Zachos. 
\newblock {\em Nambu dynamics, deformation quantization, and superintegrability.} 
\newblock In {\it Superintegrability in classical and quantum systems,} 29-46,
\newblock CRM Proc. Lecture Notes, 37, Amer. Math. Soc., Providence, RI, 2004. 


\bibitem[DSZ]{debellis:10}J. DeBellis, C. S\"amann, R. Szabo.
\newblock {\em Quantized Nambu-Poisson manifolds and n-Lie algebras.} 
\newblock J. Math. Phys. 51 (2010), no. 12, 122303, 34 pp. 

\bibitem[DFST]{dito:97}G. Dito, M. Flato, D. Sternheimer, L. Takhtajan. 
\newblock {\em Deformation quantization and Nambu mechanics.} 
\newblock Comm. Math. Phys. 183 (1997), no. 1, 1-22. 

\bibitem[F]{filippov:85}V. Filippov.
\newblock {\em n-Lie algebras.} 
\newblock Siberian Math. J. 26 (1985), no. 6, 879-891. 

\bibitem[FU]{foth:07}T. Foth, A. Uribe. 
\newblock {\em The manifold of compatible almost complex structures and geometric quantization.} 
\newblock Comm. Math. Phys. 274 (2007), no. 2, 357-379. 

\bibitem[G]{gautheron:96}P. Gautheron. 
\newblock {\em Some remarks concerning Nambu mechanics.} 
\newblock Lett. Math. Phys. 37 (1996), no. 1, 103-116. 

\bibitem[He]{helein:04}F. H\'elein. 
\newblock {\em Hamiltonian formalisms for multidimensional calculus of variations and perturbation theory.}  
\newblock In {\it Noncompact problems at the intersection of geometry, analysis, and topology,} 
\newblock 127-147, Contemp. Math., 350, Amer. Math. Soc., Providence, RI, 2004. 

\bibitem[Hu]{huybr:10}D. Huybrechts. 
\newblock {\em Hyperk\"ahler manifolds and sheaves.} 
\newblock Proceedings of the International Congress of Mathematicians. 
\newblock Volume II, 450-460, Hindustan Book Agency, New Delhi, 2010. 

\bibitem[ILMM]{ibanez:97}R. Ib\`a\~nez, M. de Le\'on, J. Marrero, D. Mart\'in de Diego. 
\newblock {\em Dynamics of generalized Poisson and Nambu-Poisson brackets.} 
\newblock J. Math. Phys. 38 (1997), no. 5, 2332-2344. 
 
\bibitem[KV]{kaledin:99}M. Verbitsky, D. Kaledin. 
\newblock {\em Hyperkahler manifolds.}
\newblock Mathematical Physics (Somerville), 12. 
\newblock International Press, Somerville, MA, 1999.

\bibitem[KS]{karabegov:01}A. Karabegov, M. Schlichenmaier. 
\newblock {\em Identification of Berezin-Toeplitz deformation quantization.} 
\newblock J. Reine Angew. Math. 540 (2001), 49-76. 

\bibitem[MM1]{ma:07}X.Ma, G. Marinescu.  
\newblock {\em Holomorphic Morse inequalities and Bergman kernels.} 
\newblock Progress in Mathematics, 254. 
\newblock Birkh\"auser Verlag, Basel, 2007. 

\bibitem[MM2]{ma:08}X.Ma, G. Marinescu.  
\newblock {\em Toeplitz operators on symplectic manifolds.}
\newblock J. Geom. Anal. 18 (2008), no. 2, 565-611. 

\bibitem[MM3]{ma:11}X.Ma, G. Marinescu. 
\newblock {\em Berezin-Toeplitz quantization and its kernel expansion.}
\newblock Geometry and quantization, 125-166,
\newblock Trav. Math., 19, Univ. Luxemb., Luxembourg, 2011. 

\bibitem[MSw]{madsen:12}T. Madsen, A. Swann. 
\newblock {\em Multi-moment maps.} 
\newblock Adv. Math. 229 (2012), no. 4, 2287-2309. 

\bibitem[M]{martin:88}G. Martin. 
\newblock {\em A Darboux theorem for multi-symplectic manifolds.}
\newblock Lett. Math. Phys. 16 (1988), no. 2, 133-138. 

\bibitem[MSu]{mukunda:76}N. Mukunda, E. Sudarshan. 
\newblock {\em Relation between Nambu and Hamiltonian mechanics.} 
\newblock Phys. Rev. D (3) 13 (1976), no. 10, 2846-2850. 

\bibitem[N]{nambu:73}Y. Nambu. 
\newblock {\em Generalized Hamiltonian dynamics.}
\newblock Phys. Rev. D (3) 7 (1973), 2405-2412. 

\bibitem[P]{polterovich:12}L. Polterovich. 
\newblock {\em Quantum unsharpness and symplectic rigidity.} 
\newblock Lett. Math. Phys. 102 (2012), no. 3, 245-264. 

\bibitem[R1]{rogers:12}C. Rogers. 
\newblock {\em $L_{\infty}$-algebras from multisymplectic geometry.}
\newblock Lett. Math. Phys. 100 (2012), no. 1, 29-50. 

\bibitem[R2]{rogers:13}C. Rogers. 
\newblock {\em 2-plectic geometry, Courant algebroids, and categorified prequantization.}
\newblock J. Symplectic Geom. 11 (2013), no. 1, 53-91.

\bibitem[RZ]{rubins:12}Y. Rubinstein, S. Zelditch. 
\newblock {\em The Cauchy problem for the homogeneous Monge-Amp\`ere equation, I. Toeplitz quantization.} 
\newblock J. Differential Geom. 90 (2012), no. 2, 303-327. 

\bibitem[SS]{samann:13}C. S\"amann, R. Szabo. 
\newblock {\em Groupoids, loop spaces and quantization of 2-plectic manifolds.} 
\newblock Reviews in Mathematical Physics Vol. 25 (2013), No. 03, 1330005

\bibitem[S1]{schlich:00}M. Schlichenmaier. 
\newblock {\em Deformation quantization of compact K\"ahler manifolds by Berezin-Toeplitz quantization.} 
\newblock Conf\'erence Mosh\'e Flato 1999, Vol. II (Dijon), 289-306,
\newblock Math. Phys. Stud., 22, Kluwer Acad. Publ., Dordrecht, 2000. 


\bibitem[S2]{schlich:10}M. Schlichenmaier. 
\newblock {\em Berezin-Toeplitz quantization for compact K\"ahler 
manifolds. A review of results.} 
\newblock Adv. Math. Phys., 2010, Article ID 927280, 38 pages.

\bibitem[T]{takhtajan:94}L. Takhtajan. 
\newblock {\em On foundation of the generalized Nambu mechanics.} 
\newblock Comm. Math. Phys. 160 (1994), no. 2, 295-315. 

\bibitem[V]{vaisman:99}I. Vaisman. 
\newblock {\em A survey on Nambu-Poisson brackets.}
\newblock Acta Math. Univ. Comenian. (N.S.) 68 (1999), no. 2, 213-241.

\end{thebibliography}
\end{document}